\documentclass{amsart}
\usepackage{fullpage} %To fit more on a page without having to set margins manually.
\usepackage{microtype}
\usepackage{graphicx} 
\usepackage{amsmath,dsfont}
\usepackage{amsthm}
\usepackage{amssymb}
\usepackage{xcolor} 
\usepackage{colortbl} 
\usepackage{tikz} 
\usepackage{hf-tikz} 
\usepackage{makecell}
\usepackage{mathrsfs}
\usepackage{graphicx}
\usepackage{enumerate}
\usepackage{chngcntr}

\usepackage{tikz}
\usepackage[colorlinks,
linkcolor=red,
anchorcolor=blue,
citecolor=green
]{hyperref}
\newtheorem{theorem}{Theorem}
\newtheorem*{theorem*}{Theorem}
\numberwithin{theorem}{section}
\numberwithin{equation}{section}
\newtheorem{definition}[theorem]{Definition}
\newtheorem{proposition}[theorem]{Proposition}
\newtheorem{example}{Example}[section]

\newtheorem{lemma}[theorem]{Lemma}

\numberwithin{prop}{section}
\counterwithout{table}{section}

\newcommand{\R}{\mathbb R}
\newcommand{\C}{\mathbb C}
\newcommand{\Ftwo}{\mathbb F_2}

\newcommand{\sgn}{\operatorname{sgn}}

\DeclareMathOperator{\lab}{lab}

\newcommand{\vN}{\operatorname{vN}}

\newcommand{\F}{\mathcal F}
\newcommand{\cH}{\mathcal H}

\title{Factoriality of Mixed $q$-Gaussian Algebras at the Boundary}
\author{Hang Yang}
\address{School of Mathematical Sciences, Zhejiang University, Hangzhou 310027, PR China}
\email{yanghangmath@zju.edu.cn}

\begin{document}
\raggedbottom
\begin{abstract}
For arbitrary index set $I$, we determine the center and characterize factoriality for mixed $q$-Gaussian von Neumann algebras $\Gamma_Q(\cH_\R)$ associated with real symmetric matrices $Q=(q_{ij})_{i,j\in I}$ satisfying $|q_{ij}|\le1$ through the strict graph we introduced in this article. Two vertices $i$ and $j$ of $I$ are joined by an edge in the strict graph if and only if $|q_{ij}|<1$. Let $S$ denote its isolated vertices. We prove that $\Gamma_Q(\cH_\R)$ is a factor if and only if $q_{ii}=-1$ for every $i\in S$, no two-vertex component $\{a,b\}$ satisfies $q_{aa}=q_{bb}=-1$, and no nonempty finite subset $A\subseteq S$ satisfies $\prod_{j\in A\setminus \{k\}}q_{jk}=1$ for every $k\in I$. We also determine the laws of some central elements.
\end{abstract}
\maketitle

\section{Introduction}
The paper studies the center of the mixed $q$-Gaussian von Neumann algebra $\mathcal M=\Gamma_Q(\cH_\R)$ with $-1\le q_{ij}\le1$. 
In the early 1990s, motivated by mathematical physics, Bo\.{z}ejko and Speicher introduced the von Neumann algebra $\Gamma_q(\cH_\R)$ generated by $q$-Gaussian variables~\cite{BS91}.
Since then, the $q$-Gaussian algebras have been studied from both probabilistic and operator algebraic viewpoints.

Speicher~\cite{Spe92} established a noncommutative central limit theorem in which $q$-Gaussian distributions arise as limit laws. This result provides a probabilistic basis for matrix models of $q$-Gaussian algebras.

In the mixed $q$-Gaussian algebra setting, the scalar deformation parameter is replaced by a real symmetric matrix $Q=(q_{ij})$ with $|q_{ij}|\le 1$. 
Such relations arise in
Speicher's work on generalized statistics \cite{Spe93}. Bo\.{z}ejko and Speicher~\cite{BS94} constructed Fock representations for more general deformed commutation relations associated with self-adjoint contractions satisfying the braid relation. For distinct indices~$i,j$, the choices~$q_{ij}=1$ and~$q_{ij}=-1$ give commutation and anticommutation relations, respectively.

An important question concerning the associated von Neumann algebras is their factoriality.
Bo\.{z}ejko, K\"ummerer, and
Speicher \cite{BKS97} developed the Gaussian functor and proved
factoriality in infinite dimension. \'{S}niady \cite{Sni04}
obtained factoriality and there is no property Gamma under a
dimension bound depending on $q$. Ricard \cite{Ric05} proved that
$\Gamma_q(\cH_{\R})$ is a factor whenever $-1<q<1$ and
$\dim\cH_{\R}\geq2$. Skalski and Wang~\cite{SW18} proved that mixed $q$-Gaussian von Neumann algebras are factors whenever the underlying real Hilbert space has dimension at least two and the real symmetric deformation matrix $Q=(q_{ij})$ satisfies $\sup_{i,j}|q_{ij}|<1$. Their arguments use a diffuse generator algebra
and the decay obtained by moving an outside letter through a long Wick word, we call this technical approarch the \emph{local decay estimate}.

For more general commutation relations, Kr\'olak studied Wick products
and factoriality in the Yang--Baxter setting \cite{Kro00,Kro06}.
Junge and Zeng \cite{JZ24} developed an ultraproduct approach to
mixed $q$-Gaussian algebras, including the endpoint values in the
construction, and established analytic and approximation properties
under the hypotheses of their respective results. Nelson and Zeng
\cite{NZ16} proved an isomorphism with free group factors for finitely
many generators, when the matrix of coefficients is
sufficiently small. Wang
\cite{Wan19} studied the singularity of generator subalgebras under
the uniform condition $\sup_{i,j}|q_{ij}|<1$.

In this paper, we define the \emph{strict graph} $G_s$ with vertex set $I$ as follows.
\begin{definition}\label{def:Strict graph and retained components}
     Two vertices \(i,j\in I\) belong to the same connected component if and only if there exist an integer \(m\ge 0\) and vertices
$i=i_0,i_1,\dots,i_m=j$
such that
$|q_{i_{r-1}i_r}|<1, \quad r=1,\dots,m.$
    The diagonal entries are recorded separately. Let
    \begin{align*}
        S&=\{i:\{i\} \text{ is a connected component of } G_s  \},\\
        \mathcal P&=\{\{a,b\}:\{a,b\} \text{ is a component of } G_s, q_{aa}=q_{bb}=-1\},\\
        V&=S\cup \bigcup_{P\in\mathcal P} P.
    \end{align*}
\end{definition}
Our factoriality criterion, which will be proved in Theorem~\ref{thm:factor-criterion}, is as follows.
\begin{theorem}~\label{thm:factor-criterion}
    The mixed $q$-Gaussian von Neumann algebra
$\mathcal M$ is a factor if and only if the following conditions hold:
\begin{enumerate}
    \item $q_{ii}=-1$ for every $i\in S$.
    \item The strict graph $G_s$ has no two-vertex component
    $\{a,b\}$ satisfying $q_{aa}=q_{bb}=-1$.
    \item There is no nonempty finite subset $A\subseteq S$ such that
    \begin{equation*}
        \prod_{j\in A\setminus\{k\}}q_{kj}=1
        \qquad\text{for every }k\in I\,.
         \end{equation*}
\end{enumerate}
An empty product is understood to be one.
\end{theorem}

The main idea is local decay estimates and propagation along the strict graph of $Q$ reduce the support of central elements to singleton components and two-vertex components $\{a,b\}$ with $q_{aa}=q_{bb}=-1$. We then describe the center in terms of radial variables and sign monomials, which yields a factoriality criterion for mixed $q$-Gaussian von Neumann algebras.

The remainder of this paper is organized as follows. Section~\ref{sec:preliminaries} develops the mixed Fock construction, establishes the tracial and spectral facts needed at the endpoints, and constructs diffuse reference algebras. Section~\ref{sec:proof of main theorem} introduces the strict graph and retained components, proves that central elements are supported on the retained algebra, and uses radius and sign coordinates together with binary commutation data to obtain the center formula and factoriality criterion. 

\section{Preliminaries}\label{sec:preliminaries}
\subsection{The mixed Fock construction}
Let $I$ be an arbitrary index set, let $\cH_\R=\ell^2(I,\R)$, and let $\cH=\cH_\R\otimes \C$. We write $(e_i)_{i\in I}$ for its distinguished orthonormal basis.
Let $Q=(q_{ij})_{i,j\in I}$ be a real symmetric matrix with $|q_{ij}|\leq 1$ for all $i,j\in I$. Define $T_Q$ on $\cH\otimes \cH$ by
\begin{equation*}T_Q(e_i\otimes e_j)=q_{ij}e_j\otimes e_i.\end{equation*}
It is a self-adjoint contraction and satisfies the braid relation.
For $n\ge 1$, let $T_r$ act as $T_Q$ on positions $r,r+1$ and as the identity on the other positions. Let $S_n$ be the symmetric group on $n$ elements.
A reduced decomposition of $\sigma\in S_n$ defines $T_\sigma$ by the corresponding product of the $T_r$; the braid relations ensure that $T_\sigma$ is independent of the chosen reduced decomposition. 
The operator 
\begin{equation*}P_Q^{(n)}=\sum_{\sigma\in S_n} T_\sigma\end{equation*} 
is positive semidefinite by the positivity theorem in \cite{BS94}. Equip the algebraic tensor space with the inner product
\begin{equation*}\langle \xi,\eta\rangle_Q=\langle\xi,P_Q^{(n)}\eta\rangle.\end{equation*}
After quotienting its null space and completing, denote the resulting Hilbert space by $\cH_Q^{\otimes n}$. The mixed Fock space is 
\begin{equation*}\F_Q=\C \Omega \oplus\bigoplus_{n\ge 1}\cH_Q^{\otimes n}. \end{equation*}
We write $\F_0=\C\Omega$ and $\F_n=\cH_Q^{\otimes n}$ for $n\ge1$.

We write $l(e_i)$ and $l_r(e_i)$ for the left and right annihilation operators, and $l(e_i)^*$ and $l_r(e_i)^*$ for their adjoints, the creation operators. 

On algebraic tensors, the left operators are given by
\begin{align}
\begin{cases}
&l(e_i)^*\Omega=e_i,\qquad l(e_i)\Omega=0,\qquad \\
&l(e_i)^*(e_{j_1}\otimes\cdots\otimes e_{j_n})
     =e_i\otimes e_{j_1}\otimes\cdots\otimes e_{j_n},\\
 &l(e_i)(e_{j_1}\otimes\cdots\otimes e_{j_n})\\
 &\qquad=\sum_{1\le r\le n,\ r:j_r=i}\left(\prod_{s<r}q_{i j_s}\right)
       e_{j_1}\otimes\cdots\widehat{e_{j_r}}\cdots\otimes e_{j_n}.
\label{eq:ann_left}
\end{cases}
\end{align}

The right operators are given by
\begin{align}
\begin{cases}
    &l_r(e_i)^*\Omega=e_i,\qquad l_r(e_i)\Omega=0,\\
    &l_r(e_i)^*(e_{j_1}\otimes\cdots\otimes e_{j_n})
     =e_{j_1}\otimes\cdots\otimes e_{j_n}\otimes e_i,\\
&l_r(e_i)(e_{j_1}\otimes\cdots\otimes e_{j_n})\\
&\qquad=\sum_{1\le r\le n,\ r:j_r=i}\left(\prod_{s>r}q_{i j_s}\right)
       e_{j_1}\otimes\cdots\widehat{e_{j_r}}\cdots\otimes e_{j_n}.
\label{eq:ann_right}
\end{cases}
\end{align}

\begin{definition}\label{def: mixed q vna}
    The mixed $q$-Gaussian von Neumann algebra generated by the fields $W_i\equiv l(e_i)^*+l(e_i)$, $i\in I$, is denoted by $\mathcal M=\vN(W_i:i\in I)$. For $A\subseteq I$, put $\mathcal M_A=\vN(W_i:i\in A)$. We use $W^*(h)$ or $\vN(h)$ for the spectral algebra of a self-adjoint measurable operator $h$. We set $\mathcal M_\emptyset=\C1$ by convention. 
\end{definition}

Each field $W_i$ has vacuum mean zero and variance one. If it is unbounded, $\vN(W_i)$ is defined by the algebra generated by the spectral projections of its self-adjoint closure. Say a von Neumann algebra is a factor if it admits a trivial center. 

Direct substitution in (\ref{eq:ann_left}) and (\ref{eq:ann_right}) gives
\begin{align}
\begin{cases}
l(e_i)l(e_i)^*=1+q_{ii}l(e_i)^*l(e_i),\\
\Delta_i\Omega=\Omega,\qquad \Delta_i(e_{j_1}\otimes\cdots\otimes e_{j_n})=\left(\prod_{r=1}^{n}q_{i j_r}\right)e_{j_1}\otimes \cdots\otimes e_{j_n},\\
l(e_j)\Delta_i=q_{ij}\Delta_il(e_j),\qquad [l(e_i),l_r(e_j)^*]=\delta_{i,j}\Delta_i\label{equ:relation of ann and creation}.
\end{cases}
\end{align}
Here $\Delta_i$ is a bounded diagonal operator on $\F_Q$. Left and right creation operators commute, as do left and right annihilation operators. 
Consequently, the left field $W_i=l(e_i)^*+l(e_i)$ commutes on this domain with every right field $W_j^r=l_r(e_j)^*+l_r(e_j)$. 

Induction on the total degree gives
\begin{equation}\label{equ:ann_norm}
    ||l(e_i)^*|_{\F_n}||\le \sqrt{n+1},\qquad ||l(e_i)|_{\F_n}||\le \sqrt{n}.
\end{equation}
    
Indeed, the first identity in \eqref{equ:relation of ann and creation} gives
\begin{equation*}\left\|l(e_i)^*\xi\right\|^2=\langle\xi,l(e_i)l(e_i)^*\xi\rangle=\|\xi\|^2+q_{ii}\|l(e_i)\xi\|^2,\end{equation*}
which bounds the squared creation norm on degree $n$ by $1$ plus the squared creation norm on degree $n-1$; for $q_{ii}\le 0$, it is bounded by $1$. 
The same argument works on subspaces with a prescribed number of occurrences of $i$. If $\xi$ is a vector of degree at most $N$, a product of $r$ fields $W_{i_1}\cdots W_{i_r}\xi$ has norm at most 
\begin{equation}\label{bound of W1W2...Wr}
2^r\sqrt{\frac{(N+r)!}{N!}}||\xi||.
\end{equation}
Let $T$ be a possibly unbounded operator on $\cH$. A vector $\xi\in\bigcap_{n\ge0}\operatorname{Dom}(T^n)$ is called an entire analytic vector for $T$ if, for every $z\in\C$, the series $\sum_{n=0}^{\infty}\frac{z^n}{n!}T^n\xi$ converges in $\cH$.  

Every finite-degree vector is entire analytic for each individual left or right field. Indeed, for complex numbers $z_1,\ldots,z_s$, the sum of the norms of the terms in a product of $s$ field exponential series, applied to a vector $\xi$ of degree at most $N$, is bounded by~\eqref{bound of W1W2...Wr} as follows:
\begin{equation*}
\|\xi\|\sum_{m\ge0}\frac{\bigl(2\sum_{j=1}^{s}|z_j|\bigr)^m}{m!}\sqrt{\frac{(N+m)!}{N!}}<\infty\,.
\end{equation*}
Thus every field is essentially self-adjoint. For a left and a right field, absolute convergence permits interchange of their exponential series. Their commutation on finite-degree vectors therefore implies that their unitary groups commute, so their self-adjoint closures commute strongly.

For an algebraic tensor $\xi$, its Wick operator $W(\xi)$ is the polynomial determined by $W(\xi)\Omega=\xi$. Wick recursion expresses it as a finite sum of products of left creation and annihilation operators. Explicitly,
\begin{multline}
W\!\left(e_{j_1}\otimes\cdots\otimes e_{j_n}\right)
=
\sum_{A\subseteq [n]}
\left(
\prod_{\substack{r\in A^{c},\, s\in A\\ r<s}}
q_{j_rj_s}
\right)\\
l(e_{j_{r_1}})^*\cdots l(e_{j_{r_p}})^*
l(e_{j_{s_1}})\cdots l(e_{j_{s_{n-p}}}),
\end{multline}
where
\[
\begin{gathered}
\relax[n]=\{1,\ldots,n\},\qquad
A=\{r_1<\cdots<r_p\},
\\
A^{c}=\{s_1<\cdots<s_{n-p}\}.
\end{gathered}
\]
Here the positions in \(A\) are assigned to creation operators and those in
\(A^{c}\) to annihilation operators. The coefficient
\[
\prod_{\substack{r\in A^{c},\, s\in A\\ r<s}}
q_{j_rj_s}
\]
is the product of the \(q\)-parameters corresponding to the crossings
required to move all creation operators to the left of all annihilation
operators. As usual, an empty product is understood to be equal to \(1\).

The vacuum is cyclic for both the left and right algebras. Indeed, consider the polynomial vacuum vector $W_{i_1}\cdots W_{i_m}\Omega$. Since $e^{itW_i}$ is a bounded unitary in $\vN(W_i)$, the difference quotients $\frac{e^{itW_i}-1}{it}\Omega$ converge to $W_i\Omega$ as $t\to0$. 
Similarly, 
\begin{equation*}
    W_{i_1}\cdots W_{i_m}\Omega=\lim_{t_1,\cdots,t_m\to 0}\prod_{r=1}^{m}\frac{e^{it_rW_{i_r}}-I}{it_r}\Omega. 
\end{equation*}
Hence, all polynomial vacuum vectors are norm limits of finite difference quotients of products of field exponentials. The Wick vectors span a dense subspace of the Fock space. Strong commutation therefore makes the vacuum cyclic and separating for the left algebra. 
Its vector state $\tau$ defined by $\tau(x)=\langle x\Omega,\Omega\rangle$ is faithful and normal. It is also a trace. The defining moment formula is
\begin{equation}\label{equ:trace tau}
    \tau\!\left(W_{i_1}\cdots W_{i_{2n}}\right)
=
\sum_{\substack{
\pi\in \mathcal P_2(2n)\\
\text{pair indices agree} 
}}
\;
\prod_{\substack{
\{r,s\},\{u,v\}\in\pi\\
r<u<s<v
}}
q_{i_r i_u}
\end{equation}
All odd moments vanish, and $\mathcal P_2(2n)$ denotes the set of all pair partitions of $\{1,2,\ldots,2n\}$. Crossing of chords is invariant under cyclic rotation, so \eqref{equ:trace tau} is tracial on polynomials. The absolutely convergent multiple series established above permit termwise expansion of finite products of generating unitary exponentials. Thus $\tau$ is tracial on their unital $*$-algebra, and separate weak operator continuity extends this identity to $\mathcal M$.

We identify $L^2(\mathcal M,\tau)$ with the Fock space by $x\mapsto x\Omega$. For $A\subseteq I$, the trace-preserving conditional expectation $E_A:\mathcal M\to \mathcal M_A$ is its $L^2$ orthogonal projection. In Fock space this is precisely the projection onto tensors all of whose letters lie in $A$:
\begin{equation}
    E_A(W(\xi))=\begin{cases}
        W(\xi),\quad \text{the indices of $\xi$ all lie in $A$} \\
        0,\quad \text{otherwise}
    \end{cases} 
\end{equation}
Since a conditional expectation is the orthogonal projection in the trace representation, this determines $E_A$. In particular,
\begin{equation*}\cap_{\lambda}\mathcal M_{A_\lambda}=\mathcal M_{\cap_\lambda{A_{\lambda}}}.\end{equation*} 

We have shown that $(\mathcal M,\tau)$ is a finite tracial von Neumann algebra. For $1\le p<\infty$, the noncommutative space $L^p(\mathcal M,\tau)$ has norm $\|x\|_p=\tau(|x|^p)^{1/p}$; we set $L^\infty(\mathcal M,\tau)=\mathcal M$ with the operator norm. 
We interpret polynomial products of unbounded fields as measurable operators affiliated with $\mathcal M$. 
The fields have moments of every order, so H\"older's inequality places every field polynomial in every finite $L^p$. To prove density, fix $1\le p<\infty$ and an integer $k\ge1$ with $2k\ge p$. The degree estimate~\eqref{bound of W1W2...Wr} and $(kn)!\le k^{kn}(n!)^k$ give
\begin{equation*}
\|W_i^n\|_p\le\|W_i^n\|_{2k}=\|W_i^{kn}\Omega\|^{1/k}\le(2\sqrt{k})^n\sqrt{n!}\,.
\end{equation*}
For real $t$, let $P_N$ be the Taylor polynomial of $e^{itx}$ of degree $N-1$. The scalar Taylor remainder gives $\|e^{itW_i}-P_N(W_i)\|_p\le(2\sqrt{k}|t|)^N/\sqrt{N!}\to0$. Approximating each factor of a product of $s$ such unitaries in $L^{sp}$ and applying H\"older's inequality proves that their unital $*$-algebra lies in the $L^p$ closure of the field polynomials. Its $C^*$ closure is strongly dense in $\mathcal M$; Kaplansky density gives uniformly bounded approximants to every bounded element, converging in $L^2$ and hence in every finite $L^p$. Norm approximation in the $C^*$ closure and spectral truncation of arbitrary $L^p$ elements finish the proof. The same argument applies to every $\mathcal M_A$.

\subsection{Diffuse reference algebras}
\begin{lemma}\label{lem:One-variable laws}
    If $q_{ii}=-1$, then $W_i^2=1$ and its two spectral projections have trace $1/2$. If $q_{ii}>-1$, the spectral law of $W_i$ is atomless. In every case its law has no atom at zero.
\end{lemma}
\begin{proof}
    For $q_{ii}=-1$, adjacent copies of $i$ are null tensors, so $(l(e_i)^*)^2=l(e_i)^2=0$. The relations in \eqref{equ:relation of ann and creation} give $W_i^2=I$, and the mean is zero, hence both spectral projections of $W_i$ have trace $1/2$. For $q_{ii}=1$, $W_i$ has the Gaussian law described above.
    
   For $-1<q_{ii}<1$, the explicit $q$-Hermite orthogonality measure in~\cite[Theorem~1.10(1)]{BKS97} is absolutely continuous, hence atomless.

\end{proof}

The two-point case requires a different source of diffuseness. In our mixed Fock model, we show below that $c=W_aW_b+W_bW_a$ has an atomless law. In particular, $c$ has no spectral mass at $\pm2$.

\begin{lemma}\label{lamma:atomless of U,c,c^2}
    Suppose that $a\ne b$, $q_{aa}=q_{bb}=-1$, and $t=q_{ab}\in(-1,1)$. The spectral measures of $U=W_aW_b$, $c=W_aW_b+W_bW_a$, and $c^2=(W_aW_b+W_bW_a)^2$ are atomless. Moreover, the algebra generated by this pair satisfies 
    \begin{equation*}\vN(W_a,W_b)\cong L^\infty (\mu_c)\overline{\otimes}M_2(\C),\quad Z(\vN(W_a,W_b))=\vN(c).\end{equation*}
    Here $\mu_c$ is the spectral measure of $c$.
\end{lemma}
\begin{proof}
    We first prove that the unitary $U=W_aW_b$ has no spectral atoms.

    Since $q_{aa}=-1$, we claim that $e_a\otimes e_a=0$ in the mixed $q$-Fock space. 
    Indeed, on the second homogeneous tensor level, $P_Q^{(2)}=I+T$, where $T(e_i\otimes e_j)=q_{ij}e_j\otimes e_i$. Applying $T$ to $e_a\otimes e_a$ gives $T(e_a\otimes e_a)=q_{aa}e_a\otimes e_a=-e_a\otimes e_a$. Hence $P_Q^{(2)}(e_a\otimes e_a)=0$. Therefore, $\langle e_a\otimes e_a, e_a\otimes e_a\rangle_Q=0$ and $e_a\otimes e_a=0$ in the mixed $q$-Fock space. More generally, whenever two identical letters become adjacent,
\[
\cdots\otimes e_a\otimes e_a\otimes\cdots ,
\]
the corresponding tensor is also null. Therefore, every positive-degree homogeneous subspace associated with the pair $\{a,b\}$ is spanned by the two alternating tensors. For $m\ge 1$, let $\xi_m^{(a)}$ and $\xi_m^{(b)}$ be the two alternating tensors of degree $2m$ beginning with $a$ and $b$, respectively, and put $f_m=W(\xi_m^{(a)})$. Denote the reversal of $f_m$ by $f^*_m$. 
Note that $f_1=U$ and annihilation of an alternating tensor can survive only at an endpoint. Indeed, deleting an interior letter creates adjacent equal letters.  
Using (\ref{eq:ann_left}) and $W_a^2=W_b^2=1$ gives $f_1=U$ and $f_2=U^2-t$ and also the operator recurrence 
\begin{equation}\label{eq:recurrence of f_m}
    f_{m+1}=U(f_m+(-t)^mf_m^*),\quad m\ge 1.
\end{equation}
Thus $f_m$ is a Laurent polynomial in $U$. For every $|z|=1$,
\begin{equation}\label{eq:bound of f_m}
    \prod_{r=1}^{m-1}(1-|t|^r)\le |f_m(z)|\le \prod_{r=1}^{m-1}(1+|t|^r).
\end{equation}
Both bounds follow inductively from (\ref{eq:recurrence of f_m}). The upper bound is uniformly finite in $m$, and the lower bound is uniformly positive. 

Suppose that there exists a spectral atom of $U$ at $z_0\in \mathbb T$, i.e. $\mu_U(\{z_0\})=\alpha>0$. Define $g=\frac{1_{\{z_0\}}}{\sqrt{\alpha}}$ in $L^2(\mu_U)$, and define $h_m=\frac{f_m}{\|f_m\|_2}$. These functions are orthonormal because their vacuum vectors have distinct degrees. By Bessel's inequality, we have
\begin{equation*}
    \sum_{m=1}^{\infty}|\langle g,h_m\rangle|^2\le \|g\|^2=1.
\end{equation*}   
However, 
\begin{equation*}
    \langle g,h_m\rangle=\int_{\mathbb T}\frac{1_{\{z_0\}}(z)}{\sqrt{\alpha}}\frac{\overline{f_m(z)}}{\|f_m\|_2}\,d\mu_U(z)=\frac{\sqrt{\alpha}\overline{f_m(z_0)}}{\|f_m\|_2},\ \text{ for any } m\ge 1.
\end{equation*}
Since the lower bound of $|f_m(z)|$ is positive by (\ref{eq:bound of f_m}), there exists a positive constant $\delta$, such that $|\langle g,h_m\rangle|>\delta$ for all $m\ge 1$.
Thus, 
\begin{equation*}\sum_{m=1}^\infty|\langle g,h_m\rangle|^2\ge \sum_{m=1}^\infty \delta^2=\infty,\end{equation*}
which contradicts Bessel's inequality. Therefore, $U$ has no spectral atoms. The maps $z\mapsto z+\bar{z}$ and $z\mapsto (z+\bar{z})^2$ have finite fibers on the circle, so $c$ and $c^2$ have no spectral atoms either.

Now the element $c$ commutes with $W_a$ and $W_b$. Put $V_0=(W_b-cW_a/2)/\sqrt{1-c^2/4}$. The relations $W_a^2=W_b^2=1$ give
\begin{equation*}(W_b-cW_a/2)^2=1-c^2/4,\quad W_a(W_b-cW_a/2)=-(W_b-cW_a/2)W_a.\end{equation*}

Since there is no mass at $c=\pm 2$, $V_0$ is a bounded symmetry.   
Moreover, $W_a$ and $V_0$ commute with $c$, and the relations
\begin{equation*}
    W_a=W_a^*,\ V_0=V_0^*,\ W_a^2=V_0^2=1,\ W_aV_0=-V_0W_a.
\end{equation*}
generate a unital copy of $M_2(\C)$.

Conversely,
$W_b=\frac{c}{2}W_a+\sqrt{1-\frac{c^{2}}{4}}\,V_{0},$
so \(\vN(W_a,W_b)=\vN(c,W_a,V_{0})\). Since \(W_a\) and \(V_{0}\) commute with \(W^*(c)\), conjugation by either symmetry fixes \(\vN(c)\) pointwise. Using \(W_aV_{0}=-V_{0}W_a\), we obtain
\[
V_{0}W_aV_{0}=-W_a,\qquad W_aV_{0}W_a=-V_{0},\qquad W_a(W_aV_{0})W_a=-W_aV_{0}.
\]
Hence, for the trace-preserving conditional expectation \(E_{\vN(c)}\),
\[
E_{\vN(c)}(W_a)=E_{\vN(c)}(V_{0})=E_{\vN(c)}(W_aV_{0})=0.
\]
Therefore \(1,W_a,V_{0},W_aV_{0}\) are linearly independent over \(\vN(c)\), and the relations \(W_a^{2}=V_{0}^{2}=1\) and \(W_aV_{0}=-V_{0}W_a\) identify their span with \(M_{2}(\mathbb C)\). Consequently,
\begin{equation*}\vN(W_a,W_b)\cong L^\infty(\mu_c)\overline{\otimes}M_2(\C),\quad Z(\vN(W_a,W_b))=\vN(c).\end{equation*} 
\end{proof}

\section{The center and factoriality theorems}\label{sec:proof of main theorem}
 Recall the definition of the strict graph in~\ref{def:Strict graph and retained components}. We note that $V$ consists of the singleton components and the two-vertex components whose fields are both symmetries. We call these components retained components. Since each retained component is a connected component of the strict graph \(G_s\), there are no strict edges between it and its complement. Hence, if \(i\) belongs to a retained component \(C\) and \(k\notin C\), then
$|q_{ik}|=1,$
so necessarily \(q_{ik}\in\{1,-1\}\). Thus distinct retained components, as well as a retained component and the rest of the index set, interact only through boundary coefficients.

\begin{example}
If all off-diagonal coefficients are strict and $|I| \ge 2$, then the strict graph is connected.
It has no retained component unless $I$ consists of exactly two symmetry vertices. If all off-diagonal coefficients are boundary values, every vertex is a singleton component. 
\end{example}

 To reduce the analysis of the center to the retained components, we use two complementary arguments. The argument~\ref{lem:localization} restricts the support of a central element to a suitable boundary neighborhood, while the propagation argument~\ref{lemma:propagation} shows that the exclusion of one vertex extends along strict edges. Together, these mechanisms reduce the analysis of the center to the retained components.

\subsection{Local compactness and central support}
\begin{lemma}\label{lem:localization}
    For $q_{ii}>-1$, set
    \begin{equation*}J_i=\{i\}\cup \{j\neq i:|q_{ij}|=1\}.\end{equation*}
    Then 
    \begin{equation*}\vN(W_i)'\cap \mathcal M\subseteq \mathcal M_{J_i}.\end{equation*}
    For a strict pair of symmetries $a,b$, set
    \begin{equation*}J_{a,b}=\{a,b\}\cup \{k\notin \{a,b\}:|q_{ka}|=|q_{kb}|=1\},\end{equation*}
    With $c=W_aW_b+W_bW_a$, one has 
    \begin{equation*}\vN(c^2)'\cap \mathcal M\subseteq \mathcal M_{J_{a,b}}.\end{equation*}
\end{lemma}
\begin{proof}
    We prove both assertions by the same argument. Let the reference algebra $\mathcal N_0$ be $\vN(W_i)$ or $\vN(W_a,W_b)$, respectively; let $J$ be $J_i$ or $J_{a,b}$.  

    For the one-generator reference algebra \(\mathcal N_0=\mathcal M_{\{i\}}\), the \(n\)-th homogeneous subspace is one-dimensional and is spanned by
\[
e_i^{\otimes n}
=
\underbrace{e_i\otimes\cdots\otimes e_i}_{n\text{ copies}},
\qquad n\ge 1.
\]

For the two-generator reference algebra \(\mathcal N_0=\mathcal M_{\{a,b\}}\), where
$q_{aa}=q_{bb}=-1,$
  the \(n\)-th homogeneous subspace is spanned by the two alternating tensors
\[
\xi_n^{(a)}
=
e_a\otimes e_b\otimes e_a\otimes e_b\otimes\cdots,
\]
and
\[
\xi_n^{(b)}
=
e_b\otimes e_a\otimes e_b\otimes e_a\otimes\cdots,
\]
each truncated after \(n\) tensor factors. 

For Wick words $x,y$ such that $x$ contains a letter $k\notin J$, define $T_{x,y}(f)=E_{\mathcal N_0}(xfy)$ for polynomials $f\in\mathcal N_0$.  
After expanding \(x\) in left normal form and the right multiplication by \(y\) in right normal form, consider a term that survives the projection onto \(L^{2}(\mathcal N_{0})\). Any left creation operator with an index outside the reference set would leave an outside letter in the final tensor and would therefore be annihilated by the projection. Likewise, any outside right annihilation acts first on a reference tensor, which contains no such letter, and hence gives zero. Thus, if the fixed word \(x\) contains an outside letter \(k\), that occurrence can contribute to a surviving term only through a left annihilation operator \(l(e_k)\).

Suppose that a surviving term contains the left annihilation \(l(e_k)\), where \(k\) lies outside the reference indices. Moving \(l(e_k)\) to the right across the right creation operators and using
\[
l(e_k) l_r(e_j)^*=l_r(e_j)^*l(e_k)+\delta_{kj}\Delta_k,
\]
we see that if \(l(e_k)\) never meets a matching \(l_r(e_k)^*\), then it eventually acts on a reference tensor, and hence gives zero. Therefore every nonzero surviving term must contain a contraction with some \(l_r(e_k)^*\), and such a contraction contributes the diagonal factor \(\Delta_k\). Thus every surviving term contains at least one factor \(\Delta_k\).
Existing diagonal factors cause no difficulty: an annihilation operator can be moved past them using $l(e_j)\Delta_k=q_{kj}\Delta_kl(e_j)$. The number of resulting terms depends only on $x,y$, not on $n$.

For a Wick word \(x=W(\xi)\), let \(\deg x\) denote the tensor length of \(\xi\).
Since the operators coming from \(x\) and \(y\) can change each reference occupation number by at most \(L=\deg x+\deg y\), an intermediate vector of initial degree \(n\) still contains at least \(\max(n-L,0)\) copies of \(i\) in the one-generator case, and at least \(\max(\lfloor n/2\rfloor-L,0)\) copies of each of \(a,b\) in the two-generator case. Since \(\Delta_k\) acts diagonally on each multidegree, its norm on the relevant intermediate subspace is therefore bounded by
\[
|q_{ki}|^{\max(n-L,0)}\quad \text{or }\quad |q_{ka}q_{kb}|^{\max(\lfloor n/2\rfloor-L,0)}, 
\]
respectively.
The relevant decay base is strictly less than $1$. By the degree estimates in (\ref{equ:ann_norm}), all remaining creation and annihilation operators contribute at most a fixed polynomial factor in \(1+n\), while every intermediate total degree differs from $n$ by at most $L$.
Hence there exist constants \(C>0\), \(\gamma\ge0\), and \(0\le \rho<1\), depending only on \(x\) and \(y\), such that
\begin{equation*}\label{equ:norm-bound of T_x,y}
\left\|T_{x,y}\big|_{L^2(\mathcal N_0)\cap \mathcal F_n}\right\|
\le
C(1+n)^\gamma
\rho^{\max(\lfloor n/2\rfloor-L,0)}.
\end{equation*}

Since each homogeneous subspace \(L^2(\mathcal N_0)\cap \mathcal F_n\) has dimension at most \(2\), choose an orthonormal basis \(\{e_{n,j}\}_{j=1}^{d_n}\) with \(d_n\le2\). Then
\[
\begin{aligned}
\|T_{x,y}\|_{\mathrm{HS}}^2
&=
\sum_{n\ge0}\sum_{j=1}^{d_n}\|T_{x,y}e_{n,j}\|_2^2
\\&\le
2\sum_{n\ge0}
\left\|T_{x,y}\big|_{L^2(\mathcal N_0)\cap\mathcal F_n}\right\|^2\\&\le
2C^2
\sum_{n\ge0}
(1+n)^{2\gamma}
\rho^{2\max(\lfloor n/2\rfloor-L,0)}
<\infty,
\end{aligned}
\]
 Hence \(T_{x,y}\) extends uniquely to a Hilbert--Schmidt operator on \(L^2(\mathcal N_0)\).

This extension agrees with the actual conditional expectation on every bounded $f\in \mathcal N_0$, even if \(x\) and \(y\) are unbounded. 
Indeed, choose reference polynomials \(f_r\) such that \(f_r\to f\) in \(L^6\). By Hölder's inequality,
\[
\|x(f_r-f)y\|_2
\le
\|x\|_6\,\|f_r-f\|_6\,\|y\|_6
\longrightarrow 0.
\]
Since the conditional expectation is contractive on \(L^2\),
\[
\|E_{\mathcal N_0}(xf_ry)-E_{\mathcal N_0}(xfy)\|_2\longrightarrow0.
\]
On the other hand, the Hilbert--Schmidt extension of \(T_{x,y}\) is bounded on \(L^2(\mathcal N_0)\), and \(f_r\to f\) also in \(L^2\); hence
\[
T_{x,y}(f_r)\longrightarrow T_{x,y}(f)
\qquad\text{in }L^2.
\]
Since \(T_{x,y}(f_r)=E_{\mathcal N_0}(xf_ry)\) for polynomials \(f_r\), passing to the limit gives
\[
T_{x,y}(f)=E_{\mathcal N_0}(xfy).
\]
Thus the Hilbert--Schmidt extension agrees with the actual conditional expectation on every bounded \(f\in\mathcal N_0\).

Let $\mathcal A_0$ now denote the diffuse abelian algebra $\vN(W_i)$ or $\vN(c^2)$, inside $\mathcal N_0$. Choose weakly null unitaries $u_n\in\mathcal A_0$. They are weakly null in $L^2(\mathcal N_0)$ as well. Compactness of the Hilbert--Schmidt extension therefore gives
\begin{equation*}\|E_{\mathcal N_0}(xu_ny)\|_2\to 0\end{equation*}
for the Wick words just considered.  

By linearity, the estimate holds for every polynomial \(x\in\ker E_J\) and polynomial \(y\). To extend it to bounded elements, let \(x\in\ker E_J\) and \(y\in \mathcal M\) be bounded. Choose field polynomials $P_r,y_r$ converging to $x,y$ in $L^4$, and put $x_r=P_r-E_J(P_r)$. Since $E_J$ preserves field polynomials and is contractive on $L^4$, we have $E_J(x_r)=0$ and $\|x_r-x\|_4\le2\|P_r-x\|_4\to0$.
Since \(\|u_n\|_\infty=1\), Hölder's inequality and the \(L^2\)-contractivity of \(E_{\mathcal N_0}\) give
\[
\begin{aligned}
\|E_{\mathcal N_0}(xu_ny)-E_{\mathcal N_0}(x_ru_ny_r)\|_2
&\le \|(x-x_r)u_ny\|_2+\|x_ru_n(y-y_r)\|_2\\
&\le \|x-x_r\|_4\|y\|_4+\|x_r\|_4\|y-y_r\|_4 .
\end{aligned}
\]
The right-hand side is independent of \(n\) and tends to \(0\) as \(r\to\infty\). It follows that
\[
\|E_{\mathcal N_0}(xu_ny)\|_2\longrightarrow0.
\]

Finally, let $z\in \mathcal A_0'\cap \mathcal M$ and $z_0=z-E_J(z)$. The conditional expectation is an $\mathcal M_J$-bimodule map and $\mathcal A_0\subseteq \mathcal M_J$, so $z_0$ commutes with $\mathcal A_0$. Therefore,
\begin{equation*}E_{\mathcal N_0}(z_0u_nz_0^*)=u_nE_{\mathcal N_0}(z_0z_0^*).\end{equation*}
The left-hand side tends to zero in $L^2$, whereas the right-hand side has constant $L^2$ norm. Thus, $E_{\mathcal N_0}(z_0z_0^*)=0$, and faithfulness of the trace gives $z_0=0$. This proves the lemma.
\end{proof}

\begin{lemma}\label{lemma:propagation}
    If $z\in Z(\mathcal M)\cap \mathcal M_A$ and $k\notin A$, then $z$ belongs to $\mathcal M_{A\setminus\{j:|q_{kj}|<1\}}$.
\end{lemma}
\begin{proof}
    For a Wick word $W(\xi)$ involving only letters from $A$, the left-right annihilation relations give
    \begin{equation}\label{eq:28}
      E_A(W_kW(\xi)W_k)\Omega=\Delta_k\xi.
    \end{equation}
    Here, on a tensor of multidegree $(m_j)_{j\in A}$, 
    \begin{equation*}\Delta_k\xi=\left(\prod_{j\in A}q_{kj}^{m_j}\right)\xi.\end{equation*}
    By $L^6$ approximation, (\ref{eq:28}) holds for bounded elements of $\mathcal M_A$ as an $L^2$ identity. This is legitimate also when $W_k$ is unbounded, since $W_k$ has every finite moment, i.e. $W_k\in L^p(\mathcal M,\tau), 1\le p<\infty$.

    Now let $z\in \mathcal M_A\cap Z(\mathcal M)$. Since $z$ commutes with $W_k$ for $k\notin A$, 
    \begin{equation*}\Delta_k(z\Omega)=E_A(W_kzW_k)\Omega=E_A(zW_k^2)\Omega=z\Omega.\end{equation*}
    On the component of multidegree $(m_j)_{j\in A}$,  $\Delta_k$ acts by the scalar $\lambda=\prod_{j\in A}q_{kj}^{m_j}$. If any strict neighbor of $k$ occurs, its absolute value is less than $1$. Hence all multidegrees containing a strict neighbor of $k$ vanish in $z\Omega$. Therefore, $z\in \mathcal M_{A\setminus\{j:|q_{kj}|<1\}}$ as required.
\end{proof}

\begin{proposition}\label{prop:central belong to M_V}
    Every central element of $\mathcal M$ belongs to $\mathcal M_V$. 
\end{proposition}
\begin{proof}
    Let \(z\in Z(\mathcal M)\). We show that every strict component not contained in \(V\) is absent from the Fock support of \(z\).

Suppose first that a strict component \(C\) contains a vertex \(i\) with
\(q_{ii}>-1\) and \(|C|\ge 2\). Choose a strict neighbor \(j\in C\) of \(i\), so $|q_{ij}|<1$.
By Lemma~\ref{lem:localization}, we have
\[
z\in \mathcal M_{J_i},
\qquad
J_i=\{i\}\cup\{l\neq i:|q_{il}|=1\},
\]
so \(j\notin J_i\) and hence \(j\) does not occur in the Fock support of \(z\).
Applying Lemma~\ref{lemma:propagation} with the outside vertex \(j\) removes every
strict neighbor of \(j\). Repeating this argument along finite strict paths
excludes every vertex of \(C\).

Now suppose that \(C\) consists entirely of symmetry vertices and
\(|C|\ge3\). Choose three distinct vertices \(a,b,k\in C\) with
\(a-b-k\) a strict path. Applying Lemma~\ref{lem:localization} to the strict pair
\(\{a,b\}\) gives
$z\in \mathcal M_{J_{a,b}},$
and since \(k\notin J_{a,b}\), the vertex \(k\) is absent from the support of
\(z\). Lemma~\ref{lemma:propagation} then excludes the whole component \(C\).

Thus every strict component except a singleton component or a two-vertex
component consisting of two symmetries is absent from the central support.
These are precisely the components defining \(V\). Therefore, using the
intersection formula for the subalgebras \(\mathcal M_A\),
$z\in \mathcal M_V.$
Hence
$Z(\mathcal M)\subseteq \mathcal M_V.$
\end{proof}

\subsection{The remaining signed algebra and its center}
Proposition~\ref{prop:central belong to M_V} shows that the center is contained in the algebra generated by the retained components. To describe the center of $\mathcal M$ explicitly, we express $\mathcal M_V$ in terms of radial coordinates and signed symmetries.

For each singleton $i\in S$ define
\begin{equation}
r_i=|W_i|, \quad s_i=\sgn(W_i).
\end{equation}
The sign is a symmetry because the one-variable law has no atom at zero. For each $P=\{a,b\}\in \mathcal P$, choose an order $(a,b)$, and put 
\begin{align}
\begin{aligned}
    c_P&=W_aW_b+W_bW_a,\quad r_P=|c_P|,\\
    x_P&=W_a,\quad z_P=\sgn(c_P),\\
    v_P&=\frac{W_b-(c_P/2)W_a}{\sqrt{1-c_P^2/4}}.
\end{aligned}
\end{align}

The denominator is interpreted through measurable functional calculus. The quotient extends to a bounded symmetry. Lemma~\ref{lamma:atomless of U,c,c^2} proves that $c_P$ has no atoms. In particular, $(W_b-(c_P/2)W_a)^2=1-c_P^2/4$ and $W_b-(c_P/2)W_a$ commutes with $\sqrt{1-c_P^2/4}$. The measurable quotient $(W_b-(c_P/2)W_a)/\sqrt{1-c_P^2/4}$ is consequently the sign of $W_b-(c_P/2)W_a$ and is a bounded symmetry. 

For a singleton, $W_i=r_is_i$. For a retained pair, we have 
\begin{equation*}W_b=\frac{r_Pz_P}{2}x_P+\sqrt{1-\frac{r_P^2}{4}}v_P.\end{equation*}
Thus, the von Neumann algebra $\mathcal M_V$ is generated by $\{r_i,\allowbreak s_i,\allowbreak x_P,\allowbreak z_P,\allowbreak v_P,\allowbreak r_P\}$, where $i\in S$ and $P\in \mathcal P$. 
We claim: $r_i$ and $r_P$ commute with every field.  
All products in this paragraph are closed products of $\tau$-measurable operators affiliated with $\mathcal M$. For $i\ne j$ and $\varepsilon=q_{ij}\in\{-1,1\}$, the degree-two Fock identity gives
\begin{equation*}
\|W_iW_j-\varepsilon W_jW_i\|_2^2=\|e_i\otimes e_j-\varepsilon e_j\otimes e_i\|_Q^2=0\,.
\end{equation*}
To pass to spectral functions, let self-adjoint measurable operators $A,B$ satisfy $AB=\varepsilon BA$. Then $A^2B=BA^2$, so the bounded operator $h=(1+A^2)^{-1}$ commutes with $B$ and hence with $(B-i)^{-1}$. Since $\vN(h)=\vN(|A|)$, the radius $|A|$ strongly commutes with $B$. If $A$ has zero kernel, its measurable polar decomposition gives $\sgn(A)=|A|^{-1}A$ and $\sgn(A)B=\varepsilon B\sgn(A)$. Consequently bounded even functions of $A$ commute with $B$, and bounded odd ones acquire the sign $\varepsilon$.

For a singleton component, this proves that $\vN(r_i)$ commutes with every field's spectral algebra. For a retained pair $P=\{a,b\}$ and $k\notin P$, the boundary relations give $W_kc_P=(q_{ka}q_{kb})c_PW_k$, while $c_P$ commutes with $W_a,W_b$. Applying the preceding resolvent argument to $c_P$ shows that every $r_i,r_P$ is affiliated with $Z(\mathcal M)$. The signed functional-calculus rules also apply to $s_i$ and $z_P$. For $v_P=\sgn(W_b-c_PW_a/2)$, its numerator acquires the sign $q_{kb}$ across $W_k$, giving the sign relations used below.

We now encode the commutation relations among $\{s_i,x_P,v_P,z_P\}$ over the field $\Ftwo$ with two elements. For an index set $J$, write $\Ftwo^{(J)}$ for the binary vectors with finite support, and $\Ftwo^J$ for all binary functions on $J$. Let $(\mathbf f_i)_{i\in I}$ be the canonical basis of $\Ftwo^{(I)}$, so $\mathbf{f}_i(j)=\delta_{ij}$ with values in $\Ftwo$. 
Let $\mathcal T$ index the symmetries $s_i$ and the triples $(x_P,v_P,z_P)$; denote them collectively by $u_\alpha,\alpha\in \mathcal T$. Define their binary support-label map $\lab:\{u_\alpha:\alpha\in \mathcal T\}\to \Ftwo^{(I)}$ by 
\begin{align*}
    &\lab(s_i)=\mathbf{f}_i,\quad \lab(x_P)=\mathbf{f}_a,\\
    &\lab(v_P)=\mathbf{f}_b,\quad \lab(z_P)=\mathbf{f}_a+\mathbf{f}_b. 
\end{align*}
 The map $\lab$ records the gauge signs carried by a symmetry. Construct a symmetric zero-diagonal binary matrix $B$ on $\mathcal T$ as follows.

 The matrix \(B=(B_{\alpha\beta})_{\alpha,\beta\in \mathcal T}\) records the
commutation relations among the sign generators. It is defined to be symmetric
with zero diagonal. For the three signs \(x_P,v_P,z_P\) associated with a retained
pair \(P\), the only nonzero local entry is
\[
B_{x_P,v_P}=B_{v_P,x_P}=1,
\]
reflecting the relation \(x_Pv_P=-v_Px_P\), while \(z_P\) commutes with both
\(x_P\) and \(v_P\). If \(u_\alpha\) and \(u_\beta\) belong to distinct retained
components, all interactions between their labels are boundary coefficients, and
we define
\[
(-1)^{B_{\alpha\beta}}
=
\prod_{i\in \operatorname{supp}(\operatorname{lab}(u_\alpha))}
\prod_{j\in \operatorname{supp}(\operatorname{lab}(u_\beta))}
q_{ij}.
\]
Thus \(B_{\alpha\beta}=0\) exactly when \(u_\alpha\) and \(u_\beta\) commute, and
\(B_{\alpha\beta}=1\) when they anticommute. Similarly, for \(k\in I\setminus V\),
the matrix \(D=(D_{k\alpha})\) records the commutation relation between the
outside field \(W_k\) and the sign \(u_\alpha\), namely
\[
(-1)^{D_{k\alpha}}
=
\prod_{j\in \operatorname{supp}(\operatorname{lab}(u_\alpha))} q_{kj}.
\]
Equivalently,
\[
u_\alpha u_\beta
=
(-1)^{B_{\alpha\beta}}u_\beta u_\alpha,
\qquad
W_k u_\alpha
=
(-1)^{D_{k\alpha}}u_\alpha W_k,\quad k\notin V.
\]
Hence \(B\) encodes the internal binary commutation data of the retained sign
algebra, while \(D\) encodes its commutation with the fields outside \(V\).

Fix a total order on $\mathcal T$. For $\eta\in\Ftwo^{(\mathcal T)}$, let $u_\eta$ be the product of the finitely many $u_\alpha$ with $\eta_\alpha=1$, taken in increasing order; the empty product is $1$. Set
\begin{equation*}K=\{\eta\in\Ftwo^{(\mathcal T)}:B\eta=0,\ D\eta=0\}.\end{equation*}
 All vectors in $K$ have finite support, even if $I$ is infinite. We use the symbols $\eta,\theta$ for vectors in $\Ftwo^{(\mathcal T)}$ and write $\eta_\alpha,\theta_\alpha$ for their binary coordinates.
The label of this monomial is the finite sum $\lab(u_\eta)=\sum_\alpha \eta_\alpha\lab(u_\alpha)$. The equation $B\eta=0$ requires $u_\eta$ to commute with the symmetries inside $\mathcal M_V$, while $D\eta=0$ requires it to commute with all fields outside $\mathcal M_V$. 
Using the coordinate description of the retained algebra, we can state the center formula.
\begin{theorem}\label{main thm:explict center}
    Every radius $r_i,r_P$ is affiliated with the center of $\mathcal M$. If $\mathcal R$ is the abelian von Neumann algebra they generate, then 
    \begin{equation*}
        Z(\mathcal M)=\vN(\mathcal R,u_\eta:\eta\in K).
    \end{equation*}
    The radii are independent in the vacuum trace. Their laws are determined explicitly by the scalar parameters as follows:
    \begin{enumerate}
        \item The law of $r_i$ is the absolute-value image of the vacuum spectral measure of the one-variable Jacobi operator with weights
        \begin{equation*}
            \sqrt{[n+1]_{q_{ii}}}, \qquad [0]_q=0,\quad [n]_q=1+q+\cdots+q^{n-1}\ (n\ge1),
        \end{equation*}
        on the basis numbered $n=0,1,\ldots$. At $q_{ii}=-1$ this law is $\delta_1$; at $q_{ii}=1$ it is the absolute value of a standard Gaussian.
        \item 
 If \(t=q_{ab}\) for \(P=\{a,b\}\), the law of \(r_P\) is the
    absolute-value image of the vacuum spectral measure of the Jacobi operator
    with zero diagonal and successive squared weights
    \begin{equation}
        \beta_n
        =
        \bigl(1+(-t)^{n-1}\bigr)\bigl(1-(-t)^n\bigr),
        \qquad n\ge 1.
    \end{equation}
    Its spectral measure is supported on \([-2,2]\) and has no atoms.
\end{enumerate}
 
\end{theorem}

We prove the~\ref{main thm:explict center} in the following subsection together with the theorem~\ref{thm:factor-criterion}. We first to describe the vacuum laws of the retained radial variables using the standard correspondence between three-term recurrences, Jacobi operators, and spectral measures; see~\cite[Sections~2.5--2.6]{Tes00}. For a singleton component, the relevant recurrence is the usual $q$-Hermite recurrence~\cite[Definition~1.9 and Theorem~1.10(1)]{BKS97}. For a retained pair~$P=\{a,b\}$, the coefficients below agree with the Al-Salam--Chihara recurrence with parameters~$(1,-1)$ and base~$-q_{ab}$ after rescaling the spectral variable; see~\cite[Section~3.8]{KS96}.

For a singleton \(i\in S\), let \(\mathsf J_i\)
be the tridiagonal operator on \(\ell^2(\mathbb N_0)\) defined by
\[
\begin{gathered}
\mathsf J_i e_n
=
\sqrt{[n+1]_{q_{ii}}}\,e_{n+1}
+
\sqrt{[n]_{q_{ii}}}\,e_{n-1},
\\
[0]_q=0,\quad [n]_q=1+q+\cdots+q^{n-1}\ (n\ge1),
\end{gathered}
\]
with \(e_{-1}=0\). Equivalently,
\[
\mathsf J_i=
\begin{pmatrix}
0 & \sqrt{[1]_{q_{ii}}} & 0 & \cdots\\
\sqrt{[1]_{q_{ii}}} & 0 & \sqrt{[2]_{q_{ii}}} & \cdots\\
0 & \sqrt{[2]_{q_{ii}}} & 0 & \ddots\\
\vdots & \vdots & \ddots & \ddots
\end{pmatrix}.
\]
Its spectral measure with respect to \(e_0\) is the vacuum law of \(W_i\), and
the law of \(r_i=|W_i|\) is its pushforward under \(x\mapsto |x|\).

Similarly, for a retained pair \(P=\{a,b\}\), with \(t=q_{ab}\), let \(\mathsf J_P\)
be the Jacobi operator with zero diagonal and off-diagonal weights
\[
\sqrt{\beta_n},
\qquad
\beta_n=
\bigl(1+(-t)^{n-1}\bigr)\bigl(1-(-t)^n\bigr),
\qquad n\ge1.
\]
Thus
\[
\mathsf J_P=
\begin{pmatrix}
0 & \sqrt{\beta_1} & 0 & \cdots\\
\sqrt{\beta_1} & 0 & \sqrt{\beta_2} & \cdots\\
0 & \sqrt{\beta_2} & 0 & \ddots\\
\vdots & \vdots & \ddots & \ddots
\end{pmatrix}.
\]
Its spectral measure at \(e_0\) is the vacuum law of
\(c_P=W_aW_b+W_bW_a\), while the law of \(r_P=|c_P|\) is again the
absolute-value pushforward of this measure.

The proof of the theorem~\ref{thm:factor-criterion}  is organized around two reductions. First, we have shown that $Z(\mathcal M)\subseteq \mathcal M_V$. Second, a trace-orthogonal expansion in radial coefficients and sign monomials turns centrality in $\mathcal M$ into the equations $B\eta=0,D\eta=0$. The spectral lemmas then identify the radius laws and decide when the radial part of the center is nontrivial.

We check the trace and independence. For polynomials $p_C$ from finitely many distinct strict components, in a product containing one contiguous block from each component, (\ref{equ:trace tau}) gives 
\begin{equation}\label{equ:product trace}
\tau\left(\prod_C p_C\right)=\prod_C\tau(p_C).
\end{equation}
Every pairing stays within one block and no crossing involves two different blocks. Approximation in $L^p$ extends the above formula to bounded spectral functions.

For a singleton component \(i\in S\), recall that
$r_i=|W_i|,\quad s_i=\operatorname{sgn}(W_i).$

Define the sign-change unitary on the mixed Fock space by
\[
U_i\bigl(e_{j_1}\otimes\cdots\otimes e_{j_n}\bigr)
=
(-1)^{m_i}
e_{j_1}\otimes\cdots\otimes e_{j_n},
\]
where \(m_i\) denotes the number of occurrences of the index \(i\).
Since changing the sign of the basis vector \(e_i\) preserves the mixed Fock
inner product, \(U_i\) is unitary and fixes the vacuum vector. Hence
$\alpha_i(x)=U_i x U_i^*$
is a trace-preserving automorphism with
$\alpha_i(W_i)=-W_i.$
Consequently,
\[
\alpha_i(r_i)=r_i,
\qquad
\alpha_i(s_i)=-s_i.
\]
Therefore, for every bounded Borel function \(f\),
\[
\tau\bigl(f(r_i)s_i\bigr)
=
\tau\bigl(\alpha_i(f(r_i)s_i)\bigr)
=
-\tau\bigl(f(r_i)s_i\bigr),
\]
and hence
\[
\tau\bigl(f(r_i)s_i\bigr)=0.
\]
Thus \(s_i\) has conditional mean zero over the radial algebra generated by
\(r_i\).

We now treat a retained pair \(P=\{a,b\}\). Recall that
$x_P^2=v_P^2=1,
\quad
x_Pv_P=-v_Px_P,$
while \(x_P\) and \(v_P\) commute with both \(r_P\) and \(z_P\). Hence, for any bounded Borel function
\(g=g(r_P,z_P)\), invariance of the trace under conjugation gives
\[
\tau(gx_P)
=
\tau(v_Pgx_Pv_P)
=
-\tau(gx_P)=0,
\]
and similarly
\[
\tau(gv_P)=0,
\qquad
\tau(gx_Pv_P)=0.
\]
Next, the gauge automorphism
\[
W_a\mapsto -W_a,\qquad W_b\mapsto W_b
\]
sends \(c_P=W_aW_b+W_bW_a\) to \(-c_P\), and therefore fixes
$r_P=|c_P|$
while sending
$z_P=\operatorname{sgn}(c_P)$
to \(-z_P\). Since this automorphism is trace-preserving,
\[
\tau\bigl(f(r_P)z_P\bigr)=0
\]
for every bounded Borel function \(f\). Thus \(z_P\) is a uniform sign and is independent of \(r_P\).

Combining this with the preceding singleton argument and the componentwise trace
factorization in \eqref{equ:product trace}, the joint law of all radii
is the product measure \(\nu\). Moreover, the ordered sign monomials satisfy
\begin{equation}\label{equ:tau for bounded radials}
\tau\!\left(f(r)u_\eta^*u_\theta\right)
=
\delta_{\eta,\theta}\int f\,d\nu.
\end{equation}
This holds for every bounded radial function \(f\). Hence the ordered sign monomials are
orthonormal over the radial algebra.

Taking \(f=1\) in (\ref{equ:tau for bounded radials}), we obtain
\[
\tau(u_\eta^*u_\theta)=\delta_{\eta,\theta},
\]
so the ordered sign monomials form an orthonormal family in
\(L^2(\mathcal M_V,\tau)\). More generally, for bounded radial functions
\(g,h\), we have
\[
\begin{aligned}
\langle g(r)u_\eta,h(r)u_\theta\rangle_2
&=\tau\bigl(g(r)\overline{h(r)}u_\theta^*u_\eta\bigr)\\
&=\delta_{\eta,\theta}\int g\overline{h}\,d\nu.
\end{aligned}
\]
By approximation, this identity extends to all
\(g,h\in L^2(\nu)\).

Moreover, since
\(u_\alpha^2=1\) and
\[
u_\alpha u_\beta
=(-1)^{B_{\alpha\beta}}u_\beta u_\alpha,
\]
every product of signs reduces to an ordered monomial
up to a scalar sign. Together with the density of bounded
radial functions in \(L^2(\nu)\), this implies that finite
linear combinations of \(f_\eta(r)u_\eta\) are dense in
\(L^2(\mathcal M_V)\). In particular,
\[
\left\|
\sum_{\eta\in F}f_\eta(r)u_\eta
\right\|_2^2
=
\sum_{\eta\in F}\|f_\eta\|_{L^2(\nu)}^2
\]
for every finite \(F\subseteq\Ftwo^{(\mathcal T)}\).

We now identify the resulting algebraic structure.
Let \(G=\Ftwo^{(\mathcal T)}\) be the group of finitely
supported binary vectors. By the commutation relations,
the product of two ordered sign monomials satisfies
\[
u_\eta u_\theta
=
\omega_B(\eta,\theta)u_{\eta+\theta},
\]
where
\[
\omega_B(\eta,\theta)
=
(-1)^{\sum_{\alpha>\beta}
B_{\alpha\beta}\eta_\alpha\theta_\beta}.
\]
Since the exponent is bilinear over \(\mathbb F_2\),
\(\omega_B\) satisfies the \(2\)-cocycle identity
\[
\omega_B(\eta,\theta)\omega_B(\eta+\theta,\gamma)
=
\omega_B(\theta,\gamma)
\omega_B(\eta,\theta+\gamma).
\]

Define the twisted left regular representation on
\(\ell^2(G)\) by
\[
\lambda_{\omega_B}(\eta)\delta_\theta
=
\omega_B(\eta,\theta)\delta_{\eta+\theta},
\]
where \(\{\delta_\theta:\theta\in G\}\) is the canonical
orthonormal basis of \(\ell^2(G)\).
The cocycle identity implies
\[
\lambda_{\omega_B}(\eta)
\lambda_{\omega_B}(\theta)
=
\omega_B(\eta,\theta)
\lambda_{\omega_B}(\eta+\theta).
\]
Thus these operators satisfy exactly the same
multiplication relations as the ordered sign monomials.

Moreover, since
$\langle f(r)u_\eta,g(r)u_\theta\rangle_2
=
\delta_{\eta,\theta}
\int f\overline{g}\,d\nu,$
 the correspondence
\[
f(r)u_\eta
\longmapsto
f\otimes\delta_\eta
\]
extends to a unitary identification
\[
L^2(\mathcal M_V)
\cong
L^2(\nu)\otimes\ell^2(G).
\]
Since this identification also preserves multiplication
and adjoints on the dense algebraic span, it yields
the trace-preserving isomorphism
\[
\mathcal M_V
\cong
L^\infty(\nu)\,\bar\otimes\,\mathcal L_{\omega_B}(G),
\]
where
\[
\mathcal L_{\omega_B}(G)
=
\vN\bigl(\lambda_{\omega_B}(\eta):\eta\in G\bigr).
\]
When \(\mathcal T\) is finite, \(\mathcal L_{\omega_B}(G)\) is precisely
the finite-dimensional signed Clifford algebra.

\subsection{Proof of the center theorem}
\begin{proof}[Proof of Theorem~\ref{main thm:explict center}]
By the preceding commutation relations, every radius is affiliated
with \(Z(\mathcal M)\). Let \(z\in Z(\mathcal M)\). Proposition~\ref{prop:central belong to M_V} implies that
\(z\in \mathcal M_V\). By the orthogonal decomposition established above,
we may write
\begin{equation}
    z=\sum_{\eta\in\mathbb F_2^{(\mathcal T)}}f_\eta(r)u_\eta,
    \qquad f_\eta\in L^2(\nu),
    \label{eq:center-fourier-expansion}
\end{equation}
where the series converges in \(L^2(\mathcal M,\tau)\). The orthogonality
relation guarantees the uniqueness of these coefficients.

We first determine which coefficients can be nonzero.
Since \(z\) is central, it commutes with every sign generator
\(u_\alpha\). By the binary commutation relations,
$u_\alpha u_\eta=(-1)^{(B\eta)_\alpha}u_\eta u_\alpha.$
Since the radial coefficients commute with all sign generators,
we obtain
\[
u_\alpha z u_\alpha
=
\sum_{\eta}
(-1)^{(B\eta)_\alpha}f_\eta(r)u_\eta.
\]
Centrality gives \(u_\alpha z u_\alpha=z\). Comparing the
unique Fourier coefficients, we conclude that
\[
f_\eta=0
\qquad\text{whenever }B\eta\neq0.
\]

Next, fix \(k\in I\setminus V\) and consider the map
\[
\Phi_k(x)=E_V(W_kxW_k),
\qquad x\in \mathcal M_V.
\]
The Wick contraction identity, together with \(L^p\)-approximation,
shows that \(\Phi_k\) agrees on \(\mathcal M_V\) with the gauge automorphism
\[
\alpha_k(W_j)=q_{kj}W_j,
\qquad j\in V.
\]
Indeed, \(q_{kj}\in\{-1,1\}\) for every \(j\in V\), so this
sign change preserves the Fock inner product and the vacuum trace.
Consequently, \(\Phi_k\) extends to a unitary operator on
\(L^2(\mathcal M_V,\tau)\).

The gauge automorphism fixes every radius and acts on the sign
generators according to their binary labels:
\[
\Phi_k(r)=r,
\qquad
\Phi_k(u_\alpha)=(-1)^{D_{k\alpha}}u_\alpha.
\]
Therefore,
\[
\Phi_k\bigl(f(r)u_\eta\bigr)
=
(-1)^{(D\eta)_k}f(r)u_\eta.
\]
Since \(\Phi_k\) is bounded on \(L^2(\mathcal M_V,\tau)\), this identity
applies coefficientwise to the expansion
\eqref{eq:center-fourier-expansion}.

On the other hand, since \(z\in Z(\mathcal M)\cap \mathcal M_V\), we have
\[
\begin{aligned}
\Phi_k(z)
=E_V(W_kzW_k)
=E_V(zW_k^2)
=zE_V(W_k^2)
=z.
\end{aligned}
\]
Here \(E_V(W_k^2)=1\), since \(k\notin V\).
The calculation remains valid when \(W_k\) is unbounded,
because \(W_k\) has finite moments of every order and the
products are well-defined in the appropriate noncommutative
\(L^p\)-spaces.

Comparing Fourier coefficients again yields
\[
f_\eta=0
\qquad\text{whenever }(D\eta)_k=1.
\]
Since this holds for every \(k\notin V\), we conclude that
\[
f_\eta=0
\qquad\text{unless }B\eta=0\text{ and }D\eta=0.
\]
Thus only indices
\[
\eta\in K:=\{\eta\in\Ftwo^{(\mathcal T)}:B\eta=0,\ D\eta=0\}
\]
can occur in the expansion of a central element.

Conversely, every radius is affiliated with \(Z(\mathcal M)\), and for
each \(\eta\in K\), the conditions \(B\eta=0\) and \(D\eta=0\)
imply that \(u_\eta\) commutes with all retained sign generators
and all fields outside \(V\). Since the retained fields are
generated by the radii and signs, \(u_\eta\) commutes with every
field \(W_j\). Hence
$u_\eta\in Z(\mathcal M).$

Set
\[
\mathcal Z_0
=
\operatorname{vN}\bigl(\mathcal R,u_\eta:\eta\in K\bigr).
\]
We have shown that \(\mathcal Z_0\subseteq Z(\mathcal M)\).
Moreover, the \(L^2\)-space of \(\mathcal Z_0\) is the closed
linear span of
\[
\{f(r)u_\eta:\eta\in K,\ f\in L^\infty(\nu)\}.
\]
Indeed, products of the generators reduce to such coefficient
monomials, and bounded radial functions are dense in
\(L^2(\nu)\). Since the expansion of \(z\) contains only indices
in \(K\), the trace-preserving conditional expectation
\(E_{\mathcal Z_0}\), viewed as the orthogonal projection onto
\(L^2(\mathcal Z_0)\), satisfies
\[
E_{\mathcal Z_0}(z)=z.
\]
Thus \(z\in\mathcal Z_0\), proving
\[
Z(\mathcal M)=\operatorname{vN}\bigl(\mathcal R,u_\eta:\eta\in K\bigr).
\]

To determine the law of $r_P=|c_P|$, write $c=c_P$ and $t=q_{ab}$, and define the degree-$2m$ tensors
\begin{equation*}
\begin{gathered}
A_m=(e_a\otimes e_b)^{\otimes m},\qquad B_m=(e_b\otimes e_a)^{\otimes m},\\
\zeta_m=A_m+B_m\quad(m\ge1),\qquad \zeta_0=\Omega\,.
\end{gathered}
\end{equation*}
Interior annihilations produce adjacent equal letters and hence null tensors. For $m\ge2$, the endpoint terms in~\eqref{eq:ann_left} give
\begin{align*}
l(e_a)l(e_b)A_m&=t^{2m-1}A_{m-1}-(-t)^mB_{m-1},\\
l(e_b)l(e_a)A_m&=A_{m-1}+(-t)^{m-1}B_{m-1}.
\end{align*}
The mixed creation-annihilation terms vanish, and the double-creation terms contribute $A_{m+1}$. At $m=1$, we have $cA_1=A_2+(1+t)\Omega$. Interchanging $a,b$ and adding, and using $c\Omega=\zeta_1$, gives
\begin{align*}
c\zeta_0&=\zeta_1,\qquad c\zeta_1=\zeta_2+2(1+t)\zeta_0,\\
c\zeta_m&=\zeta_{m+1}+\beta_m\zeta_{m-1}\quad(m\ge2),\qquad \beta_m=(1+(-t)^{m-1})(1-(-t)^m).
\end{align*}
Set $\beta_1=2(1+t)$. The vectors $\zeta_m$ have distinct homogeneous degrees. Self-adjointness and the recurrence give $\|\zeta_m\|^2=\beta_m\|\zeta_{m-1}\|^2=\prod_{j=1}^m\beta_j>0$. Their span is $\mathbb C[c]\Omega$. With $\widetilde\zeta_m=\zeta_m/\|\zeta_m\|$, the normalized basis of the vacuum cyclic subspace satisfies
\begin{equation*}
c\widetilde\zeta_m=\sqrt{\beta_{m+1}}\,\widetilde\zeta_{m+1}+\sqrt{\beta_m}\,\widetilde\zeta_{m-1}\quad(m\ge0),\qquad \beta_0=0,\quad\widetilde\zeta_{-1}=0\,.
\end{equation*}
This identifies the vacuum law of $c$ with the stated Jacobi spectral measure, and the law of $r_P$ is its pushforward under $x\mapsto|x|$.
\end{proof}

\begin{proof}[Proof of Theorem~\ref{thm:factor-criterion}]
    By Theorem~\ref{main thm:explict center}, $Z(\mathcal M)=\operatorname{vN}\bigl(\mathcal R,u_\eta:\eta\in K\bigr)$. Thus $\mathcal M$ is a factor precisely when both the radial part and the sign part of the center are trivial.

    Let $i\in S$. If $q_{ii}=-1$, then $W_i^2=1$, so $r_i=|W_i|=1$.
    Conversely, if $q_{ii}>-1$, then 
    \begin{equation*}\tau(W_i^2)=1,\quad \tau(W_i^4)=2+q_{ii}, \end{equation*}
    and hence
    \begin{equation*}\operatorname{Var}_\tau(W_i^2)=1+q_{ii}>0.\end{equation*}
    Therefore $W_i^2=r_i^2$ is nonconstant, and so $r_i$ generates a nontrivial central spectral algebra. Thus condition~1 is necessary for factoriality.

    Likewise, if \(\mathcal P\neq\varnothing\), then every retained pair \(P\in\mathcal P\) contributes the diffuse central radius \(r_P\), by Lemma~\ref{lamma:atomless of U,c,c^2}. Hence condition~2 is also necessary.

    Assume now that conditions~1 and~2 hold. Then $V=S$, every \(W_i\), \(i\in S\), is a symmetry, and all radial variables are
    constant. Thus the center can only contain sign monomials. For a finite
subset \(A\subseteq S\), write
\[
u_A=\prod_{j\in A}W_j
\]
in the fixed order. Since every \(j\in S\) is a singleton component, we have
$q_{kj}\in\{-1,1\},
\qquad k\neq j.$

Moving \(W_k\) through \(u_A\) therefore gives
\[
W_k u_A
=
\left(
\prod_{j\in A\setminus\{k\}}q_{kj}
\right)
u_AW_k.
\]
Hence \(u_A\) is central if and only if
\[
\prod_{j\in A\setminus\{k\}}q_{kj}=1
\qquad
\text{for every }k\in I.
\]
Thus condition~3 is exactly the statement that no nonzero binary vector
belongs to the kernel \(K\).

If condition~3 holds, then \(K=\{0\}\), and since all radii are constant,
Theorem~\ref{main thm:explict center} gives
$Z(\mathcal M)=\mathbb C.$

Conversely, if condition~3 fails, there exists a nonempty finite
\(A\subseteq S\) such that
\[
\prod_{j\in A\setminus\{k\}}q_{kj}=1
\qquad
\text{for every }k\in I.
\]
Then the corresponding monomial \(u_A\) commutes with every generator
\(W_k\), and hence \(u_A\in Z(\mathcal M)\). Since \(A\neq\varnothing\), the
orthogonality relation for the sign monomials implies
\[
\tau(u_A)=0,
\]
so \(u_A\) is non-scalar. Therefore \(Z(\mathcal M)\neq\mathbb C\), and \(\mathcal M\) is not
a factor.

This proves the equivalence of conditions~1--3 with factoriality, in
both the finite- and infinite-dimensional cases.
\end{proof}

\section*{Use of Artificial Intelligence}
 The author made a reasonable usage of the GPT 6.0 SOL model. We observed that the proofs of Ricard~\cite{Ric05} and Skalski and Wang~\cite{SW18} use the estimate $\bigl\|[l(e_j)^*,l_r(e_j)]\big|_{\F_n}\bigr\|\le\rho^n$, where $\rho=\sup_{i,j}|q_{ij}|<1$, to control high-degree contributions. This prompted the author to examine which products still decay, according to the indices involved, when coefficients with $|q_{ij}|=1$ are allowed, and to treat the remaining terms separately. This led to Lemmas~\ref{lamma:atomless of U,c,c^2} and~\ref{lem:localization}, used in the main proof. The main mathematical ideas originated from human reasoning and discussion. AI suggested using the field $\mathbb{F}_2$ for an explicit description of the center and helped improve the manuscript's grammar and organization. All details were checked and verified by the author.

\end{document}